\documentclass{ap-jnmp}

\usepackage{graphicx} %graphics package - should work with eps figure
\usepackage{tikz} %tikz package for manually made graphs
\usepackage{hyperref}
\usepackage{cleveref}
\DeclareMathAlphabet{\pazocal}{OMS}{zplm}{m}{n} %pazocal letters used in the text
\usepackage{gensymb} %symbols like Celsius or degree
\title{Exact Pollard-like internal water waves}

\author{Mateusz Kluczek}

\address{School of Mathematical Sciences, Western Gateway Building\\Western Road, University College Cork\\
    Cork, Ireland\\
\email{m.kluczek@umail.ucc.ie}}

\begin{document}

\maketitle
\thispagestyle{empty}

\vphantom{\vbox{%
\begin{history}
\received{24 July 2018}
%\revised{(Day Month Year)}
\accepted{6 September 2018}
%\comby{(xxxxxxxxx)}
\end{history}
}}

%---------------ABSTRACT--------------------------%
    \begin{abstract}

    In this paper we construct a new solution which represents Pollard-like three-dimensional nonlinear geophysical internal water waves. The Pollard-like solution includes the effects of the rotation of Earth and describes the internal water wave which exists at all latitudes across Earth and propagates above the thermocline. The solution is provided in Lagrangian coordinates. In the process we derive the appropriate dispersion relation for the internal water waves in a stable stratification and discuss the particles paths. An analysis of the dispersion relation for the constructed model identifies one mode of the internal water waves.

    \end{abstract}

    \ccode{2000 Mathematics Subject Classification: 74G05, 76B15, 86A05}
%-------------------------------------------------%

%-------------------------------------------------%
    \section{Introduction}\label{sec:Introduction}
%-------------------------------------------------%
The aim of this paper is to present a new exact solution which represents a nonlinear internal water wave. The solution in this study is constructed by adapting the celebrated Pollard's solution in order to successfully describe the internal water waves. In 1970, Pollard \cite{Pollard1970} presented a surface wave solution, where he extended the remarkable Gerstner's solution \cite{Gerstner1809} by including the effects of the rotation of Earth.\\
An extensive analysis of Gerstner's solution was presented in \cite{Constantin2001b,Constantin2011,Henry2008a}. Recently, there has been a significant research activity deriving Gerstner-like solutions which model various geophysical oceanic waves including equatorially-trapped surface and internal waves \cite{Constantin2012b,Constantin2013b,Constantin2014,Henry2016b,Hsu2014a} or waves in the presence of depth invariant underlying currents \cite{Henry2013,Henry2015,Henry2016,Kluczek2016,Kluczek2017,Sanjurjo2017}. Furthermore, an instability analysis of Gerstner's solution was presented in \cite{Constantin-Germain2013}. The mathematical importance of the recently derived and analysed Gerstner-like solutions is presented in a form of a review paper in \cite{Henry2017,Ionescu-Kruse2017,Johnson2017}.\\
For rotating flows in the Pollard solution a wave experiences a very slight cross-wave tilt to the wave orbital motion associated with the planetary vorticity. Therefore, the Pollard-like solution is more suitable to describe large-scale global waters rather than Gerstner's solution; since Gerstner's solution describes the motion of a particle in the vertical plane \cite{Constantin2001b,Constantin2012b,Gerstner1809,Henry2008a}, it is more adequate for flows close to the Equator where the force alternating the particles paths vanishes and the orbits are indeed vertical. The primary novel feature of this paper is we present an exact solution representing an internal water wave. The Pollard-like internal water wave solution established in this paper describes still the circular particle orbits but now the orbits lie in a plane slightly tilted from the vertical, therefore the solution is fully three-dimensional and is essentially different to the internal water wave solutions derived for the equatorial region \cite{Constantin2013b,Constantin2014}; cf. \cite{Boyd2018} for a discussion of the oceanographical relevance of these solutions.\\
The internal water waves in a stably stratified ocean may desribe the oscillation of a thermocline \cite{Constantin-Johnson2015,Cushman-Beckers2011}. The thermocline is a sharp interface separating two horizontal layers of ocean water with constant but different densities \cite{Cushman-Beckers2011,Garrison-Ellis2016,Vallis2006}. The thermocline is a phenomenon occurring also at higher latitudes, thus it is important to emphasise the need for a solution which describes the internal water waves applicable beyond the equatorial region, as is the case in this paper. The mechanism of generation of the oscillation of the thermocline is, regrettably, outside of the scope of this paper; cf. \cite{Constantin-Johnson2015,Johnson-McPhaden-Firing,Johnson2017} for a detailed study of the thermocline and its interaction with the Equatorial Undercurrent.\\
Subsequently to the work on Gerstner-like solutions, there has been developments in the analysis of Pollard's solution for surface waves \cite{Constantin-Monismith2017,Ionescu-Kruse2015b,Ionescu-Kruse2016}. A Pollard-like solution for the surface waves in the presence of mean currents and rotation was derived in a recent research paper \cite{Constantin-Monismith2017}, with an instability analysis of the Pollard-like solution presented in \cite{Ionescu-Kruse2016}. Moreover, the surface wave solution is globally dynamically possible \cite{Sanjurjo2018}. Our purpose is to modify Pollard's solution to obtain a valid model describing the nonlinear internal water waves. By empirically examining the developed solution, we hope to produce a more complete understanding of the internal oceanic flows \cite{Constantin-Johnson2016b}. We build on this analysis to identify the dispersion relation for the internal waves, desribing the oscillation of the thermocline, which may be expressed as a polynomial of degree four by a suitable non-dimensional transformation. An analysis of the polynomial identifies one mode of the internal water wave that is a standard internal gravity wave modified very slightly by Earth's rotation.

%-------------------------------------------------%
    \section{The governing equations}\label{sec:The governing equations}
%-------------------------------------------------%
The flow pattern we investigate is described in a rotating frame with the origin at a point on Earth's surface. Therefore, the $(x,y,z)$ Cartesian coordinates represent the directions of the longitude, latitude and local vertical, respectively. The governing equations for the geophysical ocean waves are given by the Euler equations \cite{Cushman-Beckers2011,Vallis2006}

	\begin{equation*}\label{eq:Governing equation 1}
        \begin{cases}
            u_t+uu_x+vu_y+wu_z+2\Omega w\cos\phi-2\Omega v\sin\phi=-\frac{1}{\rho}P_x,\\
            v_t+uv_x+vv_y+wv_z+2\Omega u\sin\phi=-\frac{1}{\rho}P_y,\\
            w_t+uw_x+vw_y+ww_z-2\Omega u\cos\phi=-\frac{1}{\rho}P_z-g,
        \end{cases}
    \end{equation*}
coupled with the equation of mass conservation

    \begin{equation*}\label{eq:Mass conservation}
        \rho_t+u\rho_x+v\rho_y+w\rho_z=0,
    \end{equation*}
together with the equation for incompressibility

    \begin{equation}\label{eq:Incompressibility}
        u_x+v_y+w_z=0.
    \end{equation}
Here $t$ is time, $\phi$ represents the latitude, $g=9.81$m$s^{-2}$ is the constant gravitational acceleration at Earth's surface, $\rho$ is the water's density, and $P$ is the pressure, while $u$, $v$ and $w$ are the respective fluid velocity components. Earth is taken to be a sphere of radius $R=6371$ km, rotating with the constant rotational speed $\Omega = 7.29\times10^{-5}$rad$\cdot s^{-1}$ round the polar axis towards the east.

    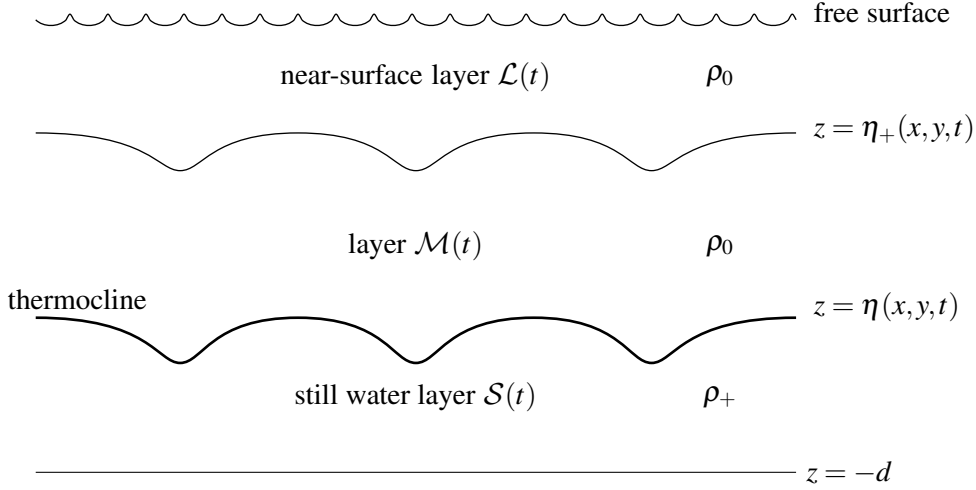
\begin{figure}[th]
        \begin{center}
            \begin{tikzpicture}

            %----------------------------------------------the free surface-----------------------------------------------------%
                \draw[line width=0.5pt] (0, 0) .. controls(0.1,-0.1)and(0.3,-0.1) .. (0.4,0) .. controls(0.45,0.1) .. (0.5, 0);
                \draw[line width=0.5pt] (0.5, 0) .. controls(0.6,-0.1)and(0.8,-0.1) .. (0.9,0) .. controls(0.95,0.1) .. (1, 0);
                \draw[line width=0.5pt] (1, 0) .. controls(1.1,-0.1)and(1.3,-0.1) .. (1.4,0) .. controls(1.45,0.1) .. (1.5, 0);
                \draw[line width=0.5pt] (1.5, 0) .. controls(1.6,-0.1)and(1.8,-0.1) .. (1.9,0) .. controls(1.95,0.1) .. (2, 0);
                \draw[line width=0.5pt] (2, 0) .. controls(2.1,-0.1)and(2.3,-0.1) .. (2.4,0) .. controls(2.45,0.1) .. (2.5, 0);
                \draw[line width=0.5pt] (2.5, 0) .. controls(2.6,-0.1)and(2.8,-0.1) .. (2.9,0) .. controls(2.95,0.1) .. (3, 0);
                \draw[line width=0.5pt] (3, 0) .. controls(3.1,-0.1)and(3.3,-0.1) .. (3.4,0) .. controls(3.45,0.1) .. (3.5, 0);
                \draw[line width=0.5pt] (3.5, 0) .. controls(3.6,-0.1)and(3.8,-0.1) .. (3.9,0) .. controls(3.95,0.1) .. (4, 0);
                \draw[line width=0.5pt] (4, 0) .. controls(4.1,-0.1)and(4.3,-0.1) .. (4.4,0) .. controls(4.45,0.1) .. (4.5, 0);
                \draw[line width=0.5pt] (4.5, 0) .. controls(4.6,-0.1)and(4.8,-0.1) .. (4.9,0) .. controls(4.95,0.1) .. (5, 0);
                \draw[line width=0.5pt] (5, 0) .. controls(5.1,-0.1)and(5.3,-0.1) .. (5.4,0) .. controls(5.45,0.1) .. (5.5, 0);
                \draw[line width=0.5pt] (5.5, 0) .. controls(5.6,-0.1)and(5.8,-0.1) .. (5.9,0) .. controls(5.95,0.1) .. (6, 0);
                \draw[line width=0.5pt] (6, 0) .. controls(6.1,-0.1)and(6.3,-0.1) .. (6.4,0) .. controls(6.45,0.1) .. (6.5, 0);
                \draw[line width=0.5pt] (6.5, 0) .. controls(6.6,-0.1)and(6.8,-0.1) .. (6.9,0) .. controls(6.95,0.1) .. (7, 0);
                \draw[line width=0.5pt] (7, 0) .. controls(7.1,-0.1)and(7.3,-0.1) .. (7.4,0) .. controls(7.45,0.1) .. (7.5, 0);
                \draw[line width=0.5pt] (7.5, 0) .. controls(7.6,-0.1)and(7.8,-0.1) .. (7.9,0) .. controls(7.95,0.1) .. (8, 0);
                \draw[line width=0.5pt] (8, 0) .. controls(8.1,-0.1)and(8.3,-0.1) .. (8.4,0) .. controls(8.45,0.1) .. (8.5, 0);
                \draw[line width=0.5pt] (8.5, 0) .. controls(8.6,-0.1)and(8.8,-0.1) .. (8.9,0) .. controls(8.95,0.1) .. (9, 0);
                \draw[line width=0.5pt] (9, 0) .. controls(9.1,-0.1)and(9.3,-0.1) .. (9.4,0) .. controls(9.45,0.1) .. (9.5, 0);
                \draw[line width=0.5pt] (9.5, 0) .. controls(9.6,-0.1)and(9.8,-0.1) .. (9.9,0) .. controls(9.95,0.1) .. (10, 0);

            %--------------------------------------------the upper boundary-----------------------------------------------------%
                \draw[line width=0.5pt] (0, -1.5) .. controls(1.5,-1.5)and(1.6,-2) .. (1.9,-2) .. controls(2.2,-2)and(2.3,-1.5) .. (3.45,-1.5) .. controls(4.6,-1.5)and(4.7,-2) .. (5,-2)  .. controls (5.3,-2)and(5.4,-1.5) .. (6.55,-1.5) .. controls(7.7,-1.5)and(7.8,-2) .. (8.1,-2)  .. controls (8.4,-2)and(8.5,-1.5) .. (10,-1.5) ;

            %---------------------------------------------the thermocline-------------------------------------------------------%
                \draw[line width=1pt] (0, -3.95) .. controls(1.5,-3.95)and(1.6,-4.55) .. (1.9,-4.55) .. controls(2.2,-4.55)and(2.3,-3.95) .. (3.45,-3.95) .. controls(4.6,-3.95)and(4.7,-4.55) .. (5,-4.55)  .. controls (5.3,-4.55)and(5.4,-3.95) .. (6.55,-3.95) .. controls(7.7,-3.95)and(7.8,-4.55) .. (8.1,-4.55)  .. controls (8.4,-4.55)and(8.5,-3.95) .. (10,-3.95) ;

            %----------------------------------------------the bed--------------------------------------------------------------%
                \draw (0,-6) -- (10,-6);

            %---------------------------------------------the names of layers---------------------------------------------------%
                \node[draw=none] at (5,-0.75) {near-surface layer $\pazocal{L}(t)$};
                \node[draw=none] at (5,-3) {layer $\pazocal{M}(t)$};
                \node[draw=none] at (5,-5) {still water layer $\pazocal{S}(t)$};

            %--------------------------------------------the density in the layers ---------------------------------------------%
                \node[draw=none] at (9,-0.75) {$\rho_0$};
                \node[draw=none] at (9,-3) {$\rho_0$};
                \node[draw=none] at (9,-5) {$\rho_+$};

            %---------------------------------------------the names of the interfaces-------------------------------------------%
                \node[draw=none,anchor=west] at (10.1,0.1) {free surface};
                \node[draw=none,anchor=west] at (10.1,-1.45) {$z=\eta_{+}(x,y,t)$};
                \node[draw=none,anchor=west] at (10.1,-3.8) {$z=\eta(x,y,t)$};
                \node[draw=none,anchor=west] at (10,-6) {$z=-d$};
                \node[draw=none,anchor=west] at (-0.5,-3.7) {thermocline};

            \end{tikzpicture}
            \caption{The depiction of the main flow regions at a fixed latitude $y$. The thermocline is described by a trochoid propagating at a speed $c$. The thermocline separates two layers of water of constant however different densities $\rho_0<\rho_+$ in a stable stratification. In the solution that we present the amplitude of the internal waves decays exponentially above the thermocline and is reduced to less then 4\% of its thermocline value at the hight of half a wave-length above the thermocline.}
            \label{fig:regions}
        \end{center}
    \end{figure}
The solution that we construct models the internal water waves describing the oscillation of a thermocline and the hydrostatic model is presented as follow. The thermocline separates layers of ocean water of different densities \cite{Cushman-Beckers2011}. The layer of less dense water $\pazocal{M}(t)$ with density $\rho_0$ overlays the layer of more dense water $\pazocal{S}(t)$ with density $\rho_+>\rho_0$. The wave motion in $\pazocal{M}(t)$ is describing the oscillations of the thermocline. The layer $\pazocal{M}(t)$ is bounded by the thermocline $z=\eta(x,y,t)$ and by the upper boundary $z=\eta_+(x,y,t)$. In the solution which we present below the amplitude of the internal waves decays exponentially with the hight above the thermocline. The amplitude of the internal waves is reduced to less then 4\% of its thermocline value at the hight of half a wave-length above the thermocline, since $e^{-\pi}\approx0.04$ (cf. \cite{Constantin2012b}), therefore, for the purposes of this model, it is justifiable to consider that the layer $\pazocal{M}(t)$ is finite and bounded. The motion in the near surface layer $\pazocal{L}(t)$ is neglected as it is a small perturbation of the free surface caused primarily by the wind and the geophysical effect has little bearing there. The layer $\pazocal{S}(t)$ of water under the thermocline describes a motionless abyssal deep-water region. The idea is to approximate the thermal structure of an ocean in the simplest form. We investigate the internal water waves in a relatively narrow ocean strip less than a few degrees of latitude wide, and so we regard the Coriolis parameters

    \begin{equation*}\label{eq:Coriolis parameters}
        f=2\Omega\sin\phi, \qquad \hat{f}=2\Omega\cos\phi,
    \end{equation*}
as constants, where $f$ is called the Coriolis parameter and $\hat{f}$ has no traditional name but usually is called the reciprocal Coriolis parameter \cite{Cushman-Beckers2011}. The typical values of the Coriolis parameters at $45\degree$ on the Northern Hemisphere are $f=\hat{f}=10^{-4}s^{-1}$ \cite{Gill1982}. On a rotating sphere, such as Earth, the Coriolis term varies with the sine of latitude, however in the $\beta$-plane approximation the Coriolis parameter is set to vary linearly in space. Furthermore, this variation can be ignored and a value of Coriolis parameter appropriate for a particular latitude can be used in the whole domain \cite{Cushman-Beckers2011}. Thus, the Euler equations reduce in the $f$-plane approximation to

    \begin{equation}\label{eq:Governing equation 2}
        \begin{cases}
            u_t+uu_x+vu_y+wu_z+\hat{f}w-fv=-\frac{1}{\rho}P_x,\\
            v_t+uv_x+vv_y+wv_z+fu=-\frac{1}{\rho}P_y,\\
            w_t+uw_x+vw_y+ww_z-\hat{f}u=-\frac{1}{\rho}P_z-g.
        \end{cases}
    \end{equation}
Water is still under the thermocline which indicates that the velocity field is in the form

    \begin{equation*}\label{eq:still water}
      (u,v,w)=(0,0,0) \mbox { for } z<\eta(x,y,t).
    \end{equation*}
Since there is no motion in the layer $\pazocal{S}(t)$ the governing equations imply the hydrostatic pressure

    \begin{equation*}\label{eq:Hydrostatic pressure}
      P=P_0-\rho_+gz \qquad z<\eta(x,y,t).
    \end{equation*}
The governing equations for the internal water waves in the layer $\pazocal{M}(t)$ are

    \begin{equation}\label{eq:Governing equation in M(t)}
        \begin{cases}
            u_t+uu_x+vu_y+wu_z+\hat{f}w-fv=-\frac{1}{\rho_0}P_x,\\
            v_t+uv_x+vv_y+wv_z+fu=-\frac{1}{\rho_0}P_y,\\
            w_t+uw_x+vw_y+ww_z-\hat{f}u=-\frac{1}{\rho_0}P_z-g.
        \end{cases}
    \end{equation}
The appropriate boundary conditions for the internal water waves are the dynamic and kinematic conditions,

    \begin{eqnarray*}\label{eq:dynamic condition}
      P=P_0-\rho_+gz \mbox{ on the thermocline } z=\eta(x,y,t)\\
	  w=\eta_t+u\eta_x+v\eta_y \mbox{ on the thermocline } z=\eta(x,y,t),
    \end{eqnarray*}
respectively. The kinematic condition prevents mixing of particles between the abyssal water region and the layer $\pazocal{M}(t)$. The particle initially on the boundary stays on the boundary all the times.

%-------------------------------------------------%
    \section{Discussion of the model}\label{sec:Discussion of the model}
%-------------------------------------------------%

%-------------------------------------------------%
        \subsection{Exact and explicit solution}\label{subsec:Exact and explicit solution}
%-------------------------------------------------%

In this section we present an exact solution to the governing equations for the internal water waves in the layer $\pazocal{M}(t)$. The Pollard-like solution represents a periodic travelling wave in the longitudinal direction at a speed of propagation $c$. For the explicit description of this flow it is convenient to use the Lagrangian framework \cite{Bennett2006}. The Lagrangian positions $(x,y,z)$ of a fluid particle are given as functions of the labelling variables $(q,r,s)$, time $t$ and real parameters $a,b,c,d,k,m$. We show that the explicit solution to the governing equations~(\ref{eq:Governing equation in M(t)}) satisfying the incompressibility condition is given by

    \begin{equation}\label{eq:Pollard explicit solution}
        \left\{
            \begin{array}{l}
                x=q-be^{-ms}\sin[k(q-ct)],\\
                y=r-de^{-ms}\cos[k(q-ct)],\\
                z=s-ae^{-ms}\cos[k(q-ct)].
             \end{array}
        \right.
    \end{equation}
The constant $k=2\pi/L$ is the wavenumber corresponding to the wavelength $L$. The parameter $q$ covers the real line, while $r\in[-r_0,r_0]$, for some $r_0$, because the solution is set up around a fixed latitude $\phi$. For every fixed value of $r\in[-r_0,r_0]$, we require $s\in[s_{0},s_{+}]$, where the choice $s=s_{0}\geq s^*>0$ represents the thermocline $z=\eta(x,y,t)$ at the latitude $\phi$, while $s=s_{+}>s_0$ prescribes the interface $z=\eta_+(x,y,t)$ separating $\pazocal{L}(t)$ and $\pazocal{M}(t)$ at the same latitude. We set the parameter of the amplitude $a>0$, wavenumber $k>0$ and for waves with amplitude decreasing above the thermocline we require $m>0$. The parameter $d$ varies from $d>0$ in the Southern Hemisphere, $d<0$ in the Northern Hemisphere to $d=0$ on the Equator since it is related to the Coriolis parameter $f$, which we show later on. Moreover, the parameters $b,c,d$ must be suitably chosen in terms of $k,m,a$.

    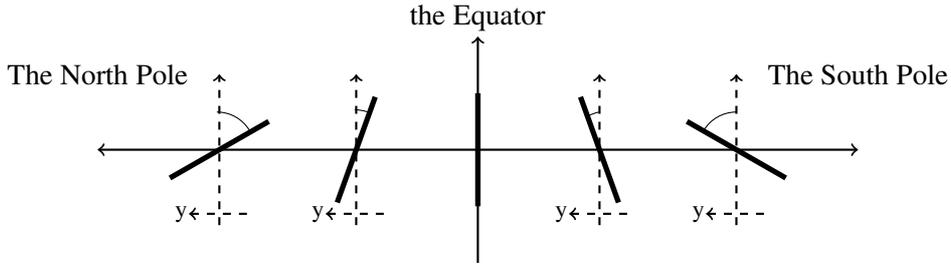
\begin{figure}[th]
        \begin{center}
            \begin{tikzpicture}

            %----------------------------------------------the names of the regions---------------------------------------------%
                \node[draw=none] at (5,-0.75) {the Equator};
                \node[draw=none] at (0,-1.5) {The North Pole};
                \node[draw=none] at (10,-1.5) {The South Pole};

            %-------------------------------------------the axis----------------------------------------------------------------%
                \draw[thick,<->] (0,-2.5) -- (10,-2.5); %latitude axis
                \draw[thick,->] (5,-4) -- (5,-1);       %vertical axis

            %-----------------------------------------------the orbits of particles---------------------------------------------%
                \draw[line width=2pt] (5,-3.25) -- (5,-1.75); %orbit at the equator
                \draw[line width=2pt] (0.95,-2.875) -- (2.25,-2.125); %first orbit from the left
                \draw[line width=2pt] (3.15,-3.2) -- (3.65,-1.8); %second orbit from the left
                \draw[line width=2pt] (6.35,-1.8) -- (6.85,-3.2); %second orbit from the right
            	\draw[line width=2pt] (7.75,-2.125) -- (9.05,-2.875); %first orbit from the right

            %------------------------------------------------the angle of inclination-------------------------------------------%
            	\draw[thick,dashed,->] (1.6,-3.5) -- (1.6,-1.5); % the angle, vertical axis and latitudinal axis for the first from the left orbit
            	\draw (2,-2.25) arc (30:90:0.5) ;
            	\draw[thick,dashed,<-] (1.2,-3.35) -- (2,-3.35);
            	\node[draw=none] at (1.1,-3.35) {\footnotesize y};
            	
            	\draw[thick,dashed,->] (3.4,-3.5) -- (3.4,-1.5); % the angle, vertical axis and latitudinal axis for the second from the left orbit
            	\draw (3.55,-2) arc (70:90:0.5) ;
            	\draw[thick,dashed,<-] (3,-3.35) -- (3.8,-3.35);
            	\node[draw=none] at (2.9,-3.35) {\footnotesize y};
            	
            	\draw[thick,dashed,->] (6.6,-3.5) -- (6.6,-1.5); % the angle, vertical axis and latitudinal axis for the second from the right orbit
            	\draw (6.6,-2) arc (100:120:0.5) ;
            	\draw[thick,dashed,<-] (6.2,-3.35) -- (7,-3.35);
            	\node[draw=none] at (6.1,-3.35) {\footnotesize y};
            	
            	\draw[thick,dashed,->] (8.4,-3.5) -- (8.4,-1.5); % the angle, vertical axis and latitudinal axis for the first from the right orbit
            	\draw (8.4,-2) arc (90:150:0.5) ;
            	\draw[thick,dashed,<-] (8,-3.35) -- (8.8,-3.35);
            	\node[draw=none] at (7.9,-3.35) {\footnotesize y};
                \end{tikzpicture}
                \caption{The figure presents the inclination of the particles orbits when the latitude increases. At the Equator the orbit becomes vertical.}
                \label{fig:Inclination}
        \end{center}
    \end{figure}
\begin{figure}[th]
        \begin{center}
            \includegraphics[width=1\textwidth]{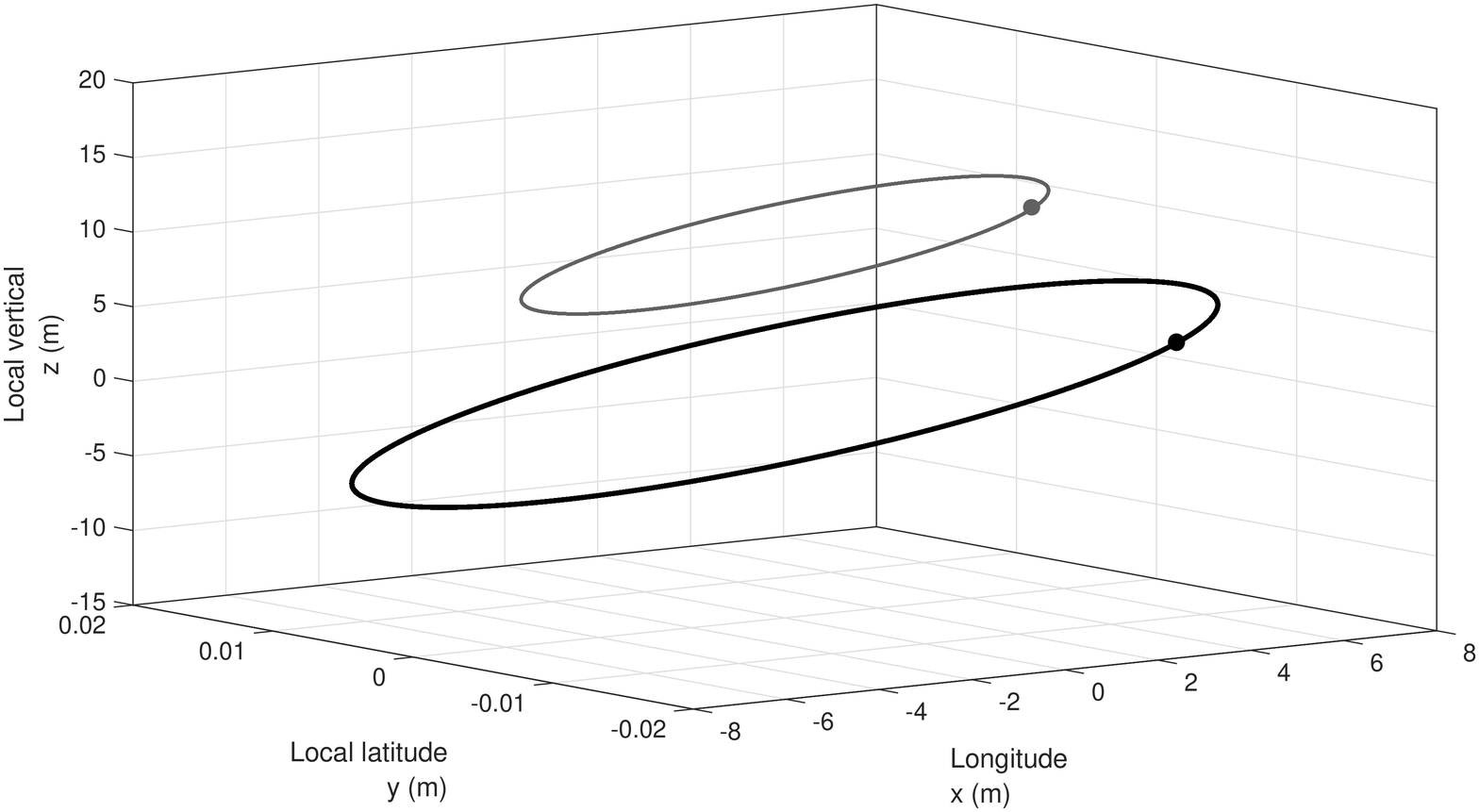}
            \caption{The path of the fluid particles when the wave propagates through water. The trajectory of the particle is a circle slightly tilted toward the Equator. The parameters of the wave-induced motion at the thermocline are $a=10$ m, $k=6.28\times 10^{-2}m^{-1}$, $\phi=45\degree$N and $\Delta\rho/\rho_0=4\times10^{-3}$. We present the motion at two depths in an ocean. The mean difference of the depths is 10 m.}
            \label{fig:Circle}
        \end{center}
    \end{figure}
\begin{figure}[th]
        \begin{center}
            \includegraphics[width=1\textwidth]{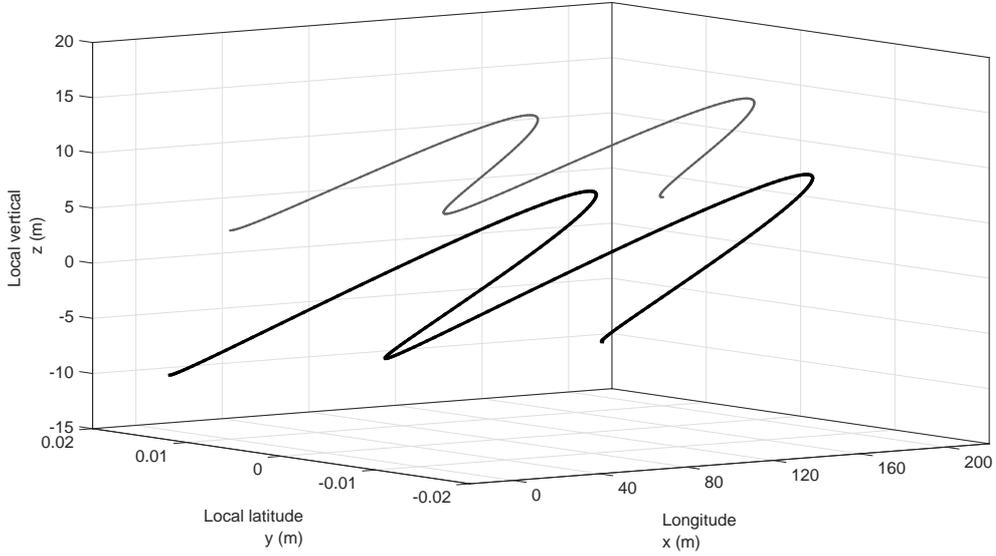}
            \caption{The trochoidal wave profile in three-dimensions generated by the oscillation of the thermocline. The wave is evaluated at two depths in an ocean (the mean difference of the depths is 10 m) for $a=10$ m, $k=6.28\times 10^{-2}m^{-1}$, $\phi=45\degree$N and $\Delta\rho/\rho_0=4\times10^{-3}$ at the thermocline. The wave profile is slightly tilted toward the Equator.}
            \label{fig:Wave_profile}
        \end{center}
    \end{figure}
\begin{figure}[th]
        \begin{center}
            \includegraphics[width=1\textwidth]{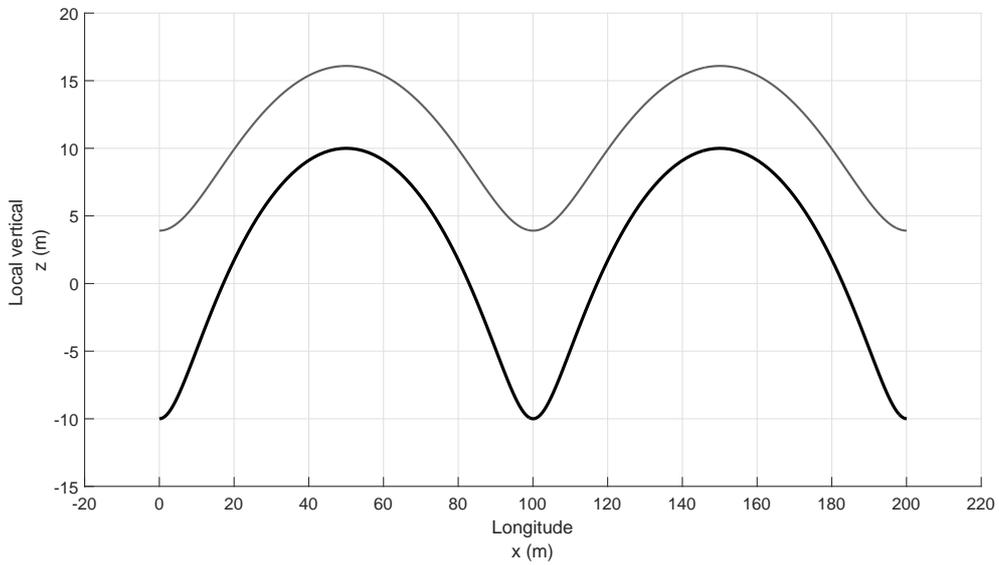}
            \caption{The projection on the vertical plane of the wave generated by the oscillation of the thermocline for two different depths (the mean difference of the depths is 10 m) at the latitude $\phi=45\degree$ on the Northern Hemisphere. The parameters of the wave at the thermocline are $a=10$ m, $k=6.28\times 10^{-2}m^{-1}$, $\Delta\rho/\rho_0=4\times10^{-3}$. The amplitude of the internal wave decreases as we ascend above the thermocline. The internal water wave is in the shape of a trochoid with narrow troughs and wide crests.}
            \label{fig:Projection}
        \end{center}
    \end{figure}

Before we proceed to proving the validity of the explicit solution~(\ref{eq:Pollard explicit solution}) we want to provide a brief discussion about the particle trajectory. For the setting of a surface wave (cf. \cite{Constantin-Monismith2017}), it is shown that the solution~(\ref{eq:Pollard explicit solution}) with parameters for the surface waves describes circles, which also applies to the internal waves. A feature of the Pollard-like solution is that the path of a particle is a slightly tilted circle \cite{Constantin-Monismith2017,Pollard1970} where the Gerstner-like solution describes circles in the vertical plane \cite{Henry2008a}. In the Pollard-like solution for the internal waves the top of the circle made by the particle is closer to the Equator and the bottom of the circle deviates to the pole at an angle of inclination $\arctan(-d/a)$ to the local vertical, which is a reversed state to the one of the surface waves \cite{Constantin-Monismith2017}. The angle of the inclination is increasing with the distance from the Equator (the figure~(\ref{fig:Inclination}). The orbits of the water particles in three-dimensions are presented in the figure~(\ref{fig:Circle}). The internal waves are in this setting in the shape of a trochoid (cf. \cite{Constantin2011}), whereas the surface wave is an inverted trochoid. The internal wave has narrow troughs and wide crests. The shape of the internal wave is depicted in the figure~(\ref{fig:Wave_profile}) taking into account the three-dimensional character. For a better explanation of the shape of the internal wave the intersection of the wave and the vertical plane is presented in the figure~(\ref{fig:Projection}). Moreover, our setting of the internal wave evaluated on the Equator particularises to the Gerstner-like equatorial internal wave solution \cite{Hsu2014a}. Note that in Gerstner's and Pollard's surface waves \cite{Constantin2001b,Constantin-Monismith2017,Henry2008a,Pollard1970} the amplitude of wave oscillations decreases as we descend into fluid, which is a reverse of the present setting, whereby the amplitude of the internal waves decreases exponentially as we ascend above the thermocline \cite{Constantin2013b,Constantin2014}. Let us now verify that~(\ref{eq:Pollard explicit solution}) is indeed the exact solution of~(\ref{eq:Governing equation in M(t)}) representing the internal water waves. For notational convenience we set

    \[
    \theta=k(q-ct).
    \]

We require

    \begin{equation*}\label{eq:Condition of r}
        s\geq s^{*}>0,
    \end{equation*}
so that $e^{-ms}<1$ throughout the layer $\pazocal{M}(t)$, since $ms\geq ks^{*}>0$. The Jacobian of the map relating the particle positions to the Lagrangian labelling variables is given by

    \begin{equation}\label{eq:the Jacobian}
        \begin{pmatrix}
            \frac{\partial x}{\partial q} & \frac{\partial y}{\partial q} & \frac{\partial z}{\partial q} \\
            \frac{\partial x}{\partial r} & \frac{\partial y}{\partial r} & \frac{\partial z}{\partial r} \\
            \frac{\partial x}{\partial s} & \frac{\partial y}{\partial s} & \frac{\partial z}{\partial s}
        \end{pmatrix}
            =
        \begin{pmatrix}
                1-kbe^{-ms}\cos\theta & kde^{-ms}\sin\theta & kae^{-ms}\sin\theta \\
                0 & 1 & 0 \\
                mbe^{-ms}\sin\theta & mde^{-ms}\cos\theta & 1+mae^{-ms}\cos\theta
        \end{pmatrix}
    \end{equation}
The flow is volume preserving and the condition of incompressibility~(\ref{eq:Incompressibility}) holds in the layer $\pazocal{M}(t)$ if and only if the determinant of the Jacobian is time independent and different than zero. The Jacobian determinant of~(\ref{eq:Pollard explicit solution}) is precisely

    \begin{equation*}\label{eq:the Jacobian determinant}
        J=1+(ma-kb)e^{-ms}\cos\theta-kmabe^{-2ms}.
    \end{equation*}
We need the condition

    \begin{equation}\label{eq:first condition}
        ma-kb=0,
    \end{equation}
to ensure that the determinant of the Jacobian is time independent. It follows that

    \begin{equation*}
        mkabe^{-2ms}\neq 1,
    \end{equation*}
throughout the flow in order to ensure a valid local diffeomorphism of~(\ref{eq:Pollard explicit solution}) by means of the inverse function theorem. Due to the condition~(\ref{eq:first condition}) and $s\geq s^{*}>0$ the above statement implies

    \begin{equation}\label{eq:ame<1}
        m^2a^2e^{-2ms^*}<1.
    \end{equation}
From the explicit solution~(\ref{eq:Pollard explicit solution}) we can deduce that the upper bound for the amplitude of internal waves is $1/m$. The Euler equations can be rewritten in the form

    \begin{equation}\label{eq:Euler eq in Lagranian}
        \begin{cases}
            \frac{Du}{Dt}+\hat{f}w-fv=-\frac{1}{\rho_{0}}P_{x},\\
            \frac{Dv}{Dt}+fu=-\frac{1}{\rho_{0}}P_{y},\\
            \frac{Dw}{Dt}-\hat{f}u=-\frac{1}{\rho_{0}}P_{z}-g,
        \end{cases}
    \end{equation}
where $D/Dt$ is the material derivative. From the direct differentiation of the system of coordinates in~(\ref{eq:Pollard explicit solution}), the velocity of each fluid particle may be expressed as

    \begin{equation}\label{eq:Velocities above thermocline}
        \left\{
            \begin{array}{l}
                u=\frac{Dx}{Dt}=kcbe^{-ms}\cos\theta,\\
                v=\frac{Dy}{Dt}=-kcde^{-ms}\sin\theta,\\
                w=\frac{Dz}{Dt}=-kcae^{-ms}\sin\theta,
            \end{array}
        \right.
    \end{equation}
and the acceleration is

    \begin{equation*}\label{eq:Acceleration above thermocline}
        \left\{
            \begin{array}{l}
                \frac{Du}{Dt}=k^2c^{2}be^{-ms}\sin\theta,\\
                \frac{Dv}{Dt}=k^2c^2de^{-ms}\cos\theta,\\
                \frac{Dw}{Dt}=k^2c^{2}ae^{-ms}\cos\theta.
            \end{array}
        \right.
    \end{equation*}
Due to the velocity and acceleration in the Lagrangian setting we can write~(\ref{eq:Euler eq in Lagranian}) as

    \begin{equation}\label{eq:Partial derivatives of pressure wrt x,y,z}
        \begin{array}{l}
            P_x=-\rho_0(k^2c^2b-kca\hat{f}+kcdf)e^{-ms}\sin\theta,\\
            P_y=-\rho_0kc(kcd+bf)e^{-ms}\cos\theta,\\
            P_z=-\rho_0(k^2c^2ae^{-ms}\cos\theta-\hat{f}kcbe^{-ms}\cos\theta+g).
        \end{array}
    \end{equation}
Since

    \begin{equation*}\label{eq:Partial derivatives of Pressure in q,s,r}
        \begin{pmatrix}
            P_q \\
            P_r \\
            P_s
        \end{pmatrix}
        =
        \begin{pmatrix}
            \frac{\partial x}{\partial q} & \frac{\partial y}{\partial q} & \frac{\partial z}{\partial q} \\
            \frac{\partial x}{\partial r} & \frac{\partial y}{\partial r} & \frac{\partial z}{\partial r} \\
            \frac{\partial x}{\partial s} & \frac{\partial y}{\partial s} & \frac{\partial z}{\partial s}
        \end{pmatrix}
        \cdot
        \begin{pmatrix}
            P_x \\
            P_y \\
            P_z
        \end{pmatrix}
    \end{equation*}
we have

    \begin{equation}\label{eq:Parital derivatives of Pressure in q,s,r 3}
        \begin{array}{l}
            P_q=-\rho_0\big[k^3c^2(a^2+d^2-b^2)e^{-ms}\cos\theta-\hat{f}kca+fkcd+k^2c^2b+kag\big]e^{-ms}\sin\theta,\\
            P_r=-\rho_0kc[kcd+bf]e^{-ms}cos\theta,\\
            P_s=-\rho_0\big[k^2c^2m(a^2+d^2-b^2)e^{-2ms}\cos^2\theta-\hat{f}kcmabe^{-2ms}+fkcmbde^{-2ms}\\
            +k^2c^2b^2me^{-2ms}+(k^2c^2a-kcb\hat{f}+mag)e^{-ms}\cos\theta+g\big].
        \end{array}
    \end{equation}
Making a natural assumption that the pressure in $\pazocal{M}(t)$ has continuous second order mixed partial derivatives we obtain the following conditions

    \begin{equation}\label{eq:second condition}
      kcd+bf=0,
    \end{equation}

    \begin{equation}\label{eq:third condition}
      mkc^2b+mcdf=k^2c^2a.
    \end{equation}
We note that the equation~(\ref{eq:second condition}) implies by means of~(\ref{eq:Partial derivatives of pressure wrt x,y,z}) that the pressure is independent of the variable $y$ throughout the layer $\pazocal{M}(t)$. Moreover, the gradient of the following pressure distribution is precisely the right-hand side of~(\ref{eq:Partial derivatives of Pressure in q,s,r})

    \begin{equation*}\label{eq:Pressure}
        \begin{array}{l}
          P=-\rho_0\big[-\frac{1}{2}k^2c^2(a^2+d^2-b^2)e^{-2ms}\cos^2\theta-\frac{1}{2}k^2c^2b^2e^{-2ms}+\frac{1}{2}\hat{f}kcabe^{-2ms}\\
          -\frac{1}{2}fkcbde^{-2ms}+(ca\hat{f}-cdf-kc^2b-ag)e^{-ms}\cos\theta+gs\big]+\tilde{P}_0.
        \end{array}
    \end{equation*}
For Pollard-like internal water waves we define that

    \begin{equation}\label{eq:P0-tildeP0}
        \begin{array}{l}
          P_0-\tilde{P}_0=-\rho_0\big[-\frac{1}{2}k^2c^2(a^2+d^2-b^2)e^{-ms_0}\cos^2\theta-\frac{1}{2}k^2c^2b^2e^{-2ms_0}\\
          +\frac{1}{2}\hat{f}kcabe^{-2ms_0}-\frac{1}{2}fkcbde^{-2ms_0}+gs_0\big]+\rho_+gs_0,
        \end{array}
    \end{equation}
to satisfy the dynamic condition. The solution $s_0$ to the equation~(\ref{eq:P0-tildeP0}) represents the thermocline. The right-hand side of~(\ref{eq:P0-tildeP0}) is a strictly increasing diffeomorphism if

    \begin{equation}\label{eq:Wavenumber}
      k>4\Omega^2/\tilde{g}\approx5\times10^{-8}\mbox{m}^{-1},
    \end{equation}
where $\tilde{g}=g(\rho_0-\rho_+)/\rho_0$ is called the coefficient of reduced gravity and $\Delta\rho/\rho_0=4\times10^{-3}$ is a typical value for the equatorial region \cite{Kessler1995}. Therefore, taking $\beta_0>P_0-\tilde{P}_0$ we can determine the solution $s_+$ representing the interface $z=\eta_+$ between the layer $\pazocal{M}(t)$ and $\pazocal{L}(t)$. Additionally, we require the continuity of pressure across the thermocline, which yields

    \begin{equation}\label{eq:continuity of pressure}
        \rho_+ga=\rho_0(kc^2b+cdf-ca\hat{f}+ag),
    \end{equation}
and pressure must be time independent hence

    \begin{equation*}\label{eq:fourth condition}
      b^2=a^2+d^2.
    \end{equation*}
From the equations~(\ref{eq:first condition}) and~(\ref{eq:second condition}) we get

    \begin{equation*}\label{eq:b}
      b=\frac{ma}{k},
    \end{equation*}

    \begin{equation*}\label{eq:d}
      d=-\frac{fma}{k^2c}.
    \end{equation*}
Therefore, the equation of continuity of the pressure~(\ref{eq:continuity of pressure}) becomes

    \begin{equation}\label{eq:continuity of pressure 2}
      \rho_0^2m^2(c^2k^2-f^2)^2=k^4(\rho_0c\hat{f}+g(\rho_+-\rho_0))^2,
    \end{equation}
Moreover, the condition~(\ref{eq:third condition}) yields

    \begin{equation*}\label{eq:m^2}
      m^2=\frac{k^4c^2}{k^2c^2-f^2},
    \end{equation*}
where $m>0$, otherwise if $m<0$ the amplitude of the wave is increasing when we ascend above the thermocline. Moreover, $m^2>0$ is ensured by~(\ref{eq:Wavenumber}) and $m=k$ at the Equator. Summarizing the aforementioned facts we obtain the dispersion relation for the internal water waves describing the oscillation of the thermocline

    \begin{equation*}\label{eq:dispersion relation}
    \rho_0^2c^2(c^2k^2-f^2)=(\rho_0c\hat{f}+g(\rho_+-\rho_0))^2.
    \end{equation*}
The dispersion relation can be simplified by including the coefficient of reduced gravity $\tilde{g}=g(\rho_0-\rho_+)/\rho_0$. Consequently, we get

    \begin{equation}\label{eq:Dispersion relation with reduced gravity}
      c^2(c^2k^2-f^2)=(c\hat{f}+\tilde{g})^2.
    \end{equation}
Choosing a suitable non-dimensional variables

    \begin{equation}\label{eq:Non-dimensional variable}
      X=c\sqrt{\frac{k}{\tilde{g}}} \qquad \varepsilon=\frac{f}{\sqrt{\tilde{g}k}} \qquad F=\frac{\hat{f}}{f},
    \end{equation}
the dispersion relation~(\ref{eq:Dispersion relation with reduced gravity}) can be rewritten as a polynomial equation of degree four $P(X)=0$ where

    \begin{equation}\label{eq:Polynomial P}
      P(X)=X^4-\varepsilon^2(1+F^2)X^2-2F\varepsilon X-1.
    \end{equation}
The roots of the polynomial $P(X)$ allows us to identify the wave speed by means of the non-dimensional variables. Moreover, we can prove that for fixed parameters there exist more than one phasespeed and we can estimate the intervals containing the roots of~(\ref{eq:Polynomial P}) (see section~(\ref{sec:Solution of the dispersion relation})). The exact value of the roots can be found by Ferrari's method. However, we focus our attention only on the existence of the real roots of the polynomial $P(X)$. The relation $c=\sqrt{\tilde{g}/k}$ refers to a standard dispersion relation for the internal waves where the Coriolis parameters are neglected \cite{Stuhlmeier2013}, which is analogous to the deep-water wave dispersion relation for surface waves \cite{Constantin2001b,Constantin2011,Constantin-Monismith2017,Henry2008a}.

%-------------------------------------------------%
        \subsection{Equatorial region}
%-------------------------------------------------%

Let us now consider the special case of a solution close to the Equator in order to substantiate the validity of the Pollard-like solution. For the equatorial waves we take the Coriolis parameters

    \begin{equation*}\label{eq:Coriolis parameters for Equatorial Waves}
      f=0, \qquad \hat{f}=2\Omega,
    \end{equation*}
and as a result, the dispersion relation~(\ref{eq:Dispersion relation with reduced gravity}) reduces to

    \begin{equation}\label{eq:Dispersion relation for equatorial waves}
      kc^2-2\Omega c-\tilde{g}=0.
    \end{equation}
The solution to the quadratic equation~(\ref{eq:Dispersion relation for equatorial waves}) is

    \begin{equation*}\label{eq:Dispersion relation for equatorial waves 2}
      c=\frac{\Omega\pm\sqrt{\Omega^2+k\tilde{g}}}{k},
    \end{equation*}
which readily agrees with the result for the internal equatorial water waves in the $f$-plane obtained in \cite{Hsu2014a}.

%-------------------------------------------------%
        \subsection{Vorticity}\label{Vorticity}
%-------------------------------------------------%

The vorticity plays important part on the trajectory of fluid particles. For an irrotational gravity-driven flow the lack of vorticity ensures that the particle paths are open loops \cite{Constantin2006,Henry2008b}. For Gerstner-like rotational flows the particle path is a closed circle \cite{Constantin2001b,Constantin2011,Henry2008a}. We prove that the Pollard-like solution we have constructed in~(\ref{eq:Pollard explicit solution}) is indeed rotational which explains somewhat the fact that the particle paths are closed circles. The vorticity is obtained by considering the product

    \begin{equation}\label{eq:Derivative of velocity wrt position variable}
      \left(\frac{\partial(q,r,s)}{\partial(x,y,z)}\right) \left(\frac{\partial(u,v,w)}{\partial(q,r,s)}\right)=\left(\frac{\partial(u,v,w)}{\partial(x,y,z)}\right),
    \end{equation}
where we exploit the inverse of~(\ref{eq:the Jacobian}) and the velocity field~(\ref{eq:Velocities above thermocline}). Moreover, the matrix~(\ref{eq:Derivative of velocity wrt position variable}) yields that the velocity field of fluid in $\pazocal{M}(t)$ is independent of the variable $y$. We are now in position to calculate the vorticity in the layer $\pazocal{M}(t)$

    \begin{equation}\label{eq:vorticity}
        \begin{array}{l}
          \omega=(w_y-v_z,u_z-w_x,v_x-u_y)=\frac{1}{1-m^2a^2e^{-2ms}}\times\\
          \begin{pmatrix}
            \frac{m^2af}{k}e^{-ms}\sin\theta\\ -c(m^2-k^2)ae^{-ms}\cos\theta+cma^2(m^2+k^2)e^{-2ms}\\ fma(\cos\theta-mae^{-ms})e^{-ms}
          \end{pmatrix}^{T}.
        \end{array}
    \end{equation}
We can validate our result by considering the vorticity in the equatorial region. Taking the Coriolis parameters for the equatorial waves $f=0$, $\hat{f}=2\Omega$, the vorticity~(\ref{eq:vorticity}) takes the form

    \begin{equation*}
        \omega=\frac{1}{1-m^2a^2e^{-2ms}}(0,2kcm^2a^2e^{-2ms},0),
    \end{equation*}
where taking the critical value of the amplitude of waves $a=1/m$ and $m=k$ we recover the vorticity for the internal equatorial water waves in the $f$-plane approximation \cite{Hsu2014a}

    \begin{equation*}
        \omega=\left(0,\frac{2kce^{-2ks}}{1-e^{-2ks}},0\right),
    \end{equation*}
and it also coincides with the vorticity in the $\beta$-plane approximation \cite{Constantin2013b}.

%-------------------------------------------------%
    \section{Solution of the dispersion relation}\label{sec:Solution of the dispersion relation}
%-------------------------------------------------%

This section presents an analytic approach towards identifying the location of roots of the polynomial~(\ref{eq:Polynomial P}). If we can find the roots of the polynomial~(\ref{eq:Polynomial P}), we can discern a wave phasespeed by means of the non-dimensional change of variables~(\ref{eq:Non-dimensional variable}). Moreover, we show that the polynomial $P(X)$, which is of degree four, has only two real roots and both are of order $O(1)$ indicating two wave speeds. It is readily seen that the constants of~(\ref{eq:Polynomial P}) are positive on both hemispheres of Earth and we can perform an analysis of the polynomial~(\ref{eq:Polynomial P}) on both hemispheres simultaneously nonetheless, we exclude the Equator since $F$ is not defined there. We recall Cauchy's theorem \cite{Prasolov2010}.

    \begin{theorem}\label{th:Cauchy theorem}
        Let $f(x)=x^n-b_1x^{n-1}-...-b_n$ where all $b_i$ are non-negative and at least one of them is non-zero. The polynomial $f$ has a unique (simple) positive root $p$ and the absolute values of the other roots do not exceed $p$.
    \end{theorem}
According to Cauchy's theorem the polynomial $P(X)$ has a unique positive root $X_0^+>0$. However, the polynomial~(\ref{eq:Polynomial P}) still can have three negative roots. We can easily compute the first derivative of the polynomial $P(X)$

    \begin{equation*}\label{eq:Derivative of polynomial P(X)}
        P'(X)=4X^3-2\varepsilon^2(1+F^2)X-2F\varepsilon
    \end{equation*}
and its discriminant

    \begin{equation*}\label{eq:Discriminant of P'(X)}
        \Delta P'(X)=128\varepsilon^6(1+F^2)^3-1728F^2\varepsilon^2.
    \end{equation*}
Making an assumption that we are outside the tropical zone, at latitudes exceeding $23\degree 26' 16"$, we have that $|F|<2.4$. Since the water temperature in the subpolar regions of Earth is constant the thermocline does not have favorable conditions to exist there and to produce the internal wave motion \cite{Garrison-Ellis2016}. Moreover, for the latitudes at most $15\degree$ away from the poles we have $|F|\geq 2-\sqrt{3}$ and therefore, we infer that the polynomial $P'(X)$ for the mid-latitudes $(23\degree 26' 16"-75\degree)$ has exactly one real root as

    \begin{equation*}\label{eq:Discriminant of P'(X) 2}
        \Delta P'(X)<0,
    \end{equation*}
which means that the polynomial $P(X)$ has one critical point. Together with $P(0)=-1$, it proves that there exist one unique positive root $X_0^+>0$ and one unique negative root $X_0^-<0$. For the polynomial $P(X)$ we can estimate

    \begin{equation}\label{eq:Estimates}
        \begin{array}{l}
            P(\pm 1)=\mp 2\varepsilon F+O(\varepsilon^2)\\
            P(1+\varepsilon F)=2\varepsilon F+O(\varepsilon^2)>0\\
            P(-1+\varepsilon F)=-2\varepsilon F+O(\varepsilon^2)<0
        \end{array}
    \end{equation}
since $F=O(1)$ and $\varepsilon=O(10^{-2})$ for internal waves with the wavelength 150-250m. Hence, the estimates~(\ref{eq:Estimates}) yield that

    \begin{equation*}
        X_0^+-1\in(0,\varepsilon F) \qquad X_0^-+1\in (0,\varepsilon F).
    \end{equation*}
for both hemispheres (see the result for the surface water waves in \cite{Constantin-Monismith2017}). We have proved therefore the existence of two real roots of the polynomial~(\ref{eq:Polynomial P}). The exact wave speed for the internal water waves generated by the oscillation of the thermocline can be found by the non-dimensional change of variable~(\ref{eq:Non-dimensional variable}) indicating two phasespeeds in dimensional terms close to

    \begin{equation*}\label{eq:Dimensional term}
      c\approx\pm\sqrt{\frac{\tilde{g}}{k}}.
    \end{equation*}
Therefore, the analysis identifies one mode of the internal wave that is a standard internal wave $c=\sqrt{\tilde{g}/k}$ \cite{Stuhlmeier2013} very slightly modified by Earth's rotation.

%-------------------------------------------------%
    \subsection*{Acknowledgements}
%-------------------------------------------------%

The author acknowledges the support of the Science Foundation Ireland (SFI) research grant 13/CDA/2117.

%---------------BIBLIOGRAPHY----------------------%

\end{document}